\documentclass%
[12pt]%
{article}

\usepackage{amssymb,amsmath, amsthm, enumerate}
\usepackage{latexsym}
\usepackage{epsfig}

\title{On Banach Spaces containing $l_p$ or $c_0$}

\author{George Androulakis \and Nigel Kalton \and Adi Tcaciuc}

\newcommand\address{\noindent\leavevmode%
\noindent Department of Mathematics,
University of South Carolina,
Columbia, SC 29208, United States\\
{\small\tt giorgis@math.sc.edu}\\[.5cm]

\noindent Department of Mathematics,
University of Missouri,
Columbia, MO 65211, United States\\
{\small\tt nigel@math.missouri.edu}\\[.5cm]

\noindent Mathematics and Statistics Department,
Grant MacEwan College,
Edmonton, Alberta, Canada T5J P2P,
{\small\tt
   tcaciuc@math.ualberta.ca}\\[.5cm]

}

\newcommand{\vphi}{\varphi}

\newcommand{\veps}{\varepsilon}

\newcommand{\de}{\delta}

\newcommand{\no}[1]{\|#1\|}

\newtheorem*{thma}{Theorem}
\newtheorem{fact}{Fact}

\newtheorem{theo}[fact]{Theorem}
\newtheorem{lem}[fact]{Lemma}

\newtheorem{coro}[fact]{Corollary}

\newcommand{\be}{\begin{equation}}
\newcommand{\ee}{\end{equation}}

\renewcommand{\thefootnote}

\usepackage{amsmath}

\begin{document}

\maketitle

\begin{abstract}
We use the Gowers block Ramsey theorem to characterize Banach spaces containing isomorphs of $\ell_p$ (for some $1 \leq p < \infty$)
or $c_0$.

\end{abstract}

\section{Introduction} \label{introduction}

 \hskip .6cm A result of Zippin \cite{Z} gives a characterization of the unit vector basis of $c_0$ and $l_p$. He showed that
a normalized basis of a Banach
space such that all normalized block bases are equivalent, must be equivalent to the unit vector basis of $c_0$ or $l_p$ for
some $1\leq p<\infty$. A Banach space $X$ with a basis $(x_i)_i$ is called \textit{asymptotic-$l_p$} \cite{MT} if there exists $K>0$ and
   an increasing function $f:\mathbb{N}\to\mathbb{N}$ such that, for all $n$, if $(y_i)_{i=1}^{n}$ is a normalized
   block basis of $(x_i)_{i=f(n)}^\infty$, then $(y_i)_{i=1}^{n}$ is equivalent to the unit vector basis of $l_p^n$.
In \cite{FFKR}  Figiel, Frankiewicz, Komorowski and Ryll-Nardzewski gave necessary and sufficient conditions for
 finding asymptotic-$l_p$ subspaces in an arbitrary Banach space. More precisely, they proved

\begin{thma}[FFKR]
 Let $p\geq 1$ and let $X$ be a Banach space with the following property:

For any infinite dimensional subspace $Y\subseteq X$ there exists a constant $M_Y$ such
that for any $n$ there exist infinite dimensional subspaces $U_1, U_2, \dots, U_n$ of $Y$ with the
property that any normalized sequence $(u_1, u_2,\dots, u_n)$ with $u_i\in U_i$ for any $i\leq n$, is
$M_Y$-equivalent to the unit vector basis of $l_p^n$.

Then $X$ contains an asymptotic-$l_p$ subspace.
\end{thma}

In \cite{Tc}
we consider similar decompositions for which any two $n$-tuples as above are uniformly
  equivalent to each other (with the equivalence constant independent of $n$) and obtain the existence of
asymptotic-$l_p$ subspaces. Our current results are in the same direction. In Theorem \ref{main} we show that if a Banach space has the property that every closed subspace contains two sequences of infinite-dimensional closed subspaces which are comparable (see the formal definition below) then it contains a copy of $\ell_p$ for some $1\leq p< \infty$ or $c_0.$  This may be regarded as a characterization of these spaces.  In Theorem \ref{main2} we show that if a Banach space $X$ is saturated with infinite dimensional subspaces $E_n$, $n \in {\mathbb N}$,
such that all normalized sequences $(x_n)$ with $x_n \in E_n$ are tail equivalent, then $X$ must contain a subspace isomorphic to $\ell_p$ or $c_0$. For the proofs we make essential use of Gowers' block Ramsey theorem \cite{G}.

 {\bf Acknowledgments.}  The third named author  wishes to thank Professor Nicole Tomczak-Jaegermann for many useful discussions.

\section{The main results}
\label{mainresult}

We first recall the Gowers block
Ramsey theorem. Let $X$ be a
Banach space with a basis $(e_i)_i$ and let $\Sigma$ be the subset
of all finite normalized block sequences of $(e_i)_i$. Given a set
$\sigma\in\Sigma$ we consider the following two players game (the Gowers
game). S (the subspace player) chooses a block subspace $X_1$ of $X$
and V (the vector player) chooses a normalized vector $x_1\in X_1$.
Then S chooses a block subspace $X_2$ and V chooses $x_2\in X_2$.
The play alternates S,V,S,V.... Then V wins the game if at some point the
sequence $(x_1,x_2\dots,x_n)$ belongs to $\sigma$.

If $Y$ is a block subspace of $X$ we say that $\sigma$ is
\emph{large} for $Y$ if every block subspace of $Y$ contains an
element of $\sigma$. We say that $\sigma$ is \emph{strategically
large} for $Y$ if player V has a winning strategy for the above game
when S is restricted to subspaces of $Y$.

Let $\Delta=(\de_1, \de_2...)$ be a sequence of positive numbers.
The set $\sigma_\Delta$, the $\Delta$-enlargement of $\sigma$, will
stand for the set of all sequences $(x_1,\dots,x_n)\in\Sigma$ such
that there exists a sequence $(y_1,\dots,y_n)\in\sigma$ with
$\no{x_i-y_i}<\de_i$, for all $i\leq n$.

\begin{theo} {\rm (Gowers's block Ramsey Theorem)} \label{Gowers}
Let $\sigma\in\Sigma$ be large for $X$ and let $\Delta$ be a
sequence of positive numbers. Then there exists a block subspace $Y$
of $X$ such that $\sigma_\Delta$ is strategically large for $Y$.
\end{theo}

The following lemma is an application of Gowers's block
Ramsey Theorem.

\begin{lem} \label{mainLemma}
Let $X$ be a Banach space with a basis. Then there exists $1\leq
p\leq\infty$ and an infinite dimensional block subspace $Y$ of $X$
such that for every sequence $(Y_i)_{i=1}^{\infty}$ of infinite
dimensional block subspaces of $Y$ and for all $n$ there exists a
normalized block sequence $(y_i)_{i=1}^{n}$ in $Y$ with $y_i\in
Y_i$ for $1\leq i\leq n$ and $(y_i)_{i=1}^{n}$ is $2-equivalent$
to the unit vector basis of $l_p^n$.

\end{lem}

For the proof of this lemma we will need to define the stabilized Krivine set of a Banach space. If $X$ is a Banach
space with a basis and $W$ is an infinite dimensional block
subspace of $X$, let $K(W)$ be the set of $p$'s in $[1,\infty]$
such that $l_p$ is finitely block represented in $W$, ($p=\infty$
if $c_0$ is finitely block represented in $W$). Then $K(W)$ is a
closed non-empty subset of
 $[1,\infty]$ \cite{K}. Note that if $W_2$ is an infinite
  dimensional block subspace of $W_1$ then $K(W_2)\subseteq K(W_1)$.
  Moreover, if $W_2$ has finite codimension in $W_1$ then $K(W_2)=K(W_1)$.
   Using these properties it can be easily shown that there exists
an infinite dimensional block $W$ of $X$ such that $K(W)=K(V)$ for
all infinite dimensional block subspaces $V$ of $W$,
 (else a transfinite induction gives a contradiction). The set $K(W)$ is called
a {\it stabilized Krivine set} for $X$.

\noindent {\bf Proof. } Let $W$ be an infinite dimensional block subspace of $X$ such that $K(W)$ is a stabilized Krivine
set for $X$. Fix $p\in K(W)$ and for any $n\in\mathbb{N}$ let $\sigma_n$ be the
set of all finite normalized block sequences of $W$ of length $n$
such that they are $2$-equivalent to the unit vector basis of
$l_p^n$. The conclusion of the Lemma follows easily if we can find
a block subspace $Y$ such that, for any $n$, $\sigma_n$ is
strategically large for $Y$.

  Note that for any $n$, $\sigma_n$
is large in any infinite dimensional block subspace $V$ of $W$. By
applying Theorem \ref{Gowers} repeatedly, we obtain a nested
sequence
 $V_1\supset V_2\supset V_3\supset\dots\supset V_n\supset\dots$ of block subspaces of $W$ such that, for any $n$, $\sigma_n$
  is strategically large for $V_n$. Note that we do not enlarge $\sigma_n$ since we can replace the $2$ in ``$2$-equivalent''
   by $1+\veps$ for any $\veps>0$.

Let $Y$ be a diagonal block subspace, that is a subspace generated
by a block basis $v_1, v_2,\dots$ with $v_n\in V_n$ for
 every $n$. We claim that $\sigma_n$ is strategically large for $Y$, for any value of $n$. Indeed fix $n\in\mathbb{N}$ and
  denote by $[Y]_n$ the $n$-tail of $Y$, that is the subspace generated by $(v_j)_{j\geq n}$. Note
that $[Y]_n\subseteq V_n$. Consider a typical Gowers game in $Y$.
For any choice of a block subspace $Z$ of $Y$ the subspace player
makes, the vector player chooses a vector $z\in Z\cap[Y]_n$ as if
the game was played inside $V_n$ and the subspace player picked
$Z\cap[Y]_n$. Since the vector player has a winning strategy for
the game played in $V_n$, it follows that after finitely many
steps the finite block sequence he chooses belongs to $\sigma_n$.
Therefore the vector player has a winning strategy for the game
played in $Y$ as well. This proves that $\sigma_n$ is
strategically large for $Y$, which finishes the proof of the
lemma.

\hfill $\blacksquare$

Let $\mathcal E=(E_j)_{j=1}^{\infty}$ be a sequence of nonzero subspaces of a Banach
space $X$  and let $\mathcal F=(F_j)_{j=1}^{\infty}$ be a sequence of nonzero subspaces
in a Banach space $Y$.  We will say that
 $\mathcal E$ is $C$-{\it dominated} by $\mathcal F$, and we write
$\mathcal E\stackrel{C}\prec\mathcal F$ if for any $n$ we have
\[
                \no{\sum_{j=1}^{n}x_j}\leq C \no{\sum_{j=1}^{n} y_j}
\]
whenever $x_j\in E_j,\ y_j\in F_j$ with $\|x_j\|=\|y_j\|$ for $1\le j\le n.$

The two sequences $\mathcal E$ and $\mathcal F$  are called $C$-{\it comparable} if
$\mathcal E\stackrel{C}{\prec}\mathcal F$ or $\mathcal F\stackrel{C}{\prec}\mathcal E$.
$\mathcal E$ and $\mathcal F$ are $C$-equivalent
if and only if there exist constants $C_1$ and $C_2$ with
$C_1C_2=C$ such that $\mathcal E\stackrel{C_1}\prec\mathcal F$
 and $\mathcal F\stackrel{C_2}{\prec}\mathcal E$. \vskip .4cm

Notice that the sequence $\mathcal E$ is comparable to itself if and only if it is
equivalent to itself and this is in turn equivalent to the fact that
$\mathcal E=(E_j)_{j=1}^\infty$ is an absolute Schauder decomposition of its
closed linear span $[\mathcal E].$

Note that $\mathcal E$ is $C$-dominated by the canonical one-dimensional decomposition of $\ell_p$ (or $c_0$ when $p=\infty$) if and only if $\mathcal E$ satisfies an upper $\ell_p$-estimate with constant $C$, i.e.:
$$ \|\sum_{j=1}^nx_j\|\le C(\sum_{j=1}^n\|x_j\|^p)^{1/p}, \qquad x_j\in E_j, \ n=1,2,\ldots$$ or
$$ \|\sum_{j=1}^nx_j\|\le C\max_{1\le j\le n}\|x_j\|, \qquad x_j\in E_j, \ n=1,2,\ldots$$ when $p=\infty$.
Similarly $\mathcal E$ $C$-dominates the canonical one-dimensional decomposition of $\ell_p$ if and only if $\mathcal E$ satisfies a lower $\ell_p-$estimate with constant $C$:
$$ (\sum_{j=1}^n\|x_j\|^p)^{1/p}\le C\|\sum_{j=1}^nx_j\|, \qquad x_j\in E_j, \ n=1,2,\ldots,$$ or
$$ \max_{1\le j\le n}\|x_j\|\le  C\|\sum_{j=1}^nx_j\|, \qquad x_j\in E_j, \ n=1,2,\ldots,$$ when $p=\infty.$

Recall that a basic sequence $(x_n)_{n=1}^{\infty}$ has an upper (respectively, lower) $p$-estimate if there is a constant $C$ so that for all block basic sequences $(u_n)_{n=1}^{\infty}$ the sequence $([u_n])_{n=1}^{\infty}$ has an upper (respectively, lower) $\ell_p$-estimate with constant $C$.

\begin{theo}\label{upperlower} Let $X$ be a separable Banach space  with a basis $(e_j)_{j=1}^{\infty}$ and suppose $1\le p\le\infty.$\newline
(i) Assume that every closed subspace of $X$ contains a sequence $\mathcal E$ of infinite-dimensional subspaces with an upper $\ell_p$-estimate.  Then $(e_j)_{j=1}^{\infty}$ has a block basic sequence with an upper $p$-estimate.\newline
(ii) Assume that every closed subspace of $X$ contains a sequence $\mathcal E$ of infinite-dimensional subspaces with a lower $\ell_p$-estimate.  Then $(e_j)_{j=1}^{\infty}$ has a block basic sequence with a lower $p$-estimate.\end{theo}

\noindent{\bf Proof.} We prove only (i) as (ii) is similar.  We may assume that every block subspace contains a sequence $\mathcal E$ of block subspaces with an upper $\ell_p-$estimate.  For each block subspace $W$ let $C(W)$ denote the infimum of all constants $C$ so that $W$ contains a sequence $\mathcal E$ of block subspaces with an upper  $\ell_p-$estimate with constant $C$. We claim that there exists an infinite dimensional block subspace
$Y$ of $X$ and a constant $C<\infty$ such that for each infinite
dimensional block subspace $Z$ of $Y$  we have that $C(Z)< C$.
Indeed, otherwise there exists a decreasing sequence of block
subspaces $Z_n$ of $X$ such that $C(Z_n)>n$. If we choose a sequence
$(w_j)_{j=1}^{\infty}$ of successive, linearly independent block vectors with
$w_j\in Z_j$ for each $j$, then the infinite dimensional block
subspace $W$ spanned by $(w_j)_{j=1}^{\infty}$ contains a sequence of infinite-dimensional subspaces $\mathcal E=(E_j)_{j=1}^{\infty}$ with an $\ell_p-$upper estimate with constant $C$, say. Picking $n>C$ and considering the sequence $(E_j\cap Z_n)_{j=1}^{\infty}$, we have a contradiction.  This contradiction proves the
above claim.

Thus we may assume that for original basis we have the property that $C(W)<C$ for every block subspace.
Now define the set $\sigma$ to consist of all normalized finite block basic sequences $(u_1,\ldots,u_n)$ so that for some $a_1,\ldots,a_n$ with $|a_1|^p+\cdots+|a_n|^p=1$
$$ \|a_1u_1+\cdots+a_nu_n\|>C+2.$$

If the conclusion of the Theorem is false, this set is large.  Let $\Delta=(2^{-i})_{i=1}^{\infty}.$  then by Theorem \ref{Gowers} $\sigma_{\Delta}$ is strategically large for some block subspace $Y$.  Let $\mathcal E=(E_j)_{j=1}^{\infty}$ be a sequence of infinite-dimensional block subspaces of $Y$ with an upper $\ell_p-$estimate with constant at most $C$.
If the subspace player S uses the strategy $\mathcal E$ then the vector player V may select normalized vectors $v_j\in E_j$ so that for some $n$ there exists $(u_1,\ldots,u_n)\in\sigma$ so that $\|u_j-v_j\|<2^{-j}$ for $j=1,2,\ldots,n.$
Now for an appropriate $a_1,\ldots,a_n$ with $|a_1|^p+\cdots+|a_n|^p=1$ we have
$$ \|\sum_{j=1}^na_jv_j\| \ge \|\sum_{j=1}^na_ju_j\|-1\ge C+1$$ which gives a contradiction.  This proves (i). \hfill $\blacksquare$

 \begin{theo} \label{main}  Let $X$ be a separable Banach space with the property
that every infinite-dimensional closed subspace contains two comparable sequences
$\mathcal E$ and $\mathcal F$ of infinite-dimensional closed subspaces.  Then $X$
contains a copy of $\ell_p$ for some $1\le p<\infty$ or $c_0.$
\end{theo}

\noindent {\bf Proof. }
We first assume that $X$ has a basis and then by Lemma \ref{mainLemma}, we may pass to a block subspace $Y$ so that for a suitable $p$, whenever $\mathcal E=(E_j)_{j=1}^\infty$ is a sequence of infinite-dimensional subspaces of $Y$ then for any $n\in\mathbb N$ there exist $y_j\in E_j$ for $j=1,2,\ldots,n$ so that $(y_j)_{j=1}^n$ is 2-equivalent to the canonical $\ell_p^n$-basis.

By our assumption, for any closed infinite dimensional subspace $Z$ of $Y$ let
${\mathcal E}$ and ${\mathcal F}$ be two comparable sequences of infinite dimensional
subspaces of $Z$ , say $\mathcal{E}= (E_j)$ is dominated by $\mathcal{F}= (F_j)$.
By Lemma~\ref{mainLemma} for any $n \in {\mathbb N}$ there exists $y_j \in E_j$ ($j=1, \ldots , n$),
such that $(y_j)_{j=1}^n$ is $2$-equivalent to the unit vector basis of $\ell_p^n$.
Thus ${\mathcal F}$ has a lower $p$ estimate. Similarly, ${\mathcal E}$ has an
upper $p$-estimate.  Applying Theorem \ref{upperlower} (i) and (ii) one after another
we can pass to a block basic sequence which has both an upper and a lower $p$-estimate, i.e. is equivalent to the canonical basis of $\ell_p.$ \hfill $\blacksquare$

\begin{coro} Let $X$ be a separable Banach space with the property that every
infinite-dimensional closed subspace contains a subspace with an absolute Schauder
decomposition of infinite-dimensional subspaces. Then $X$ contains a copy of
$\ell_p$ for some $1\le p<\infty$ or $c_0.$ \end{coro}

We remark that this Corollary could easily be deduced from the result in \cite{Tc}.

A sequence $(x_n)_n$ will be called $C$-{\it tail equivalent}
if for any $N\in\mathbb{N}$, there exists $k>N$ such that $(x_n)_{n=1}^{\infty}$
is $C$-equivalent to  $(x_n)_{n=k}^{\infty}$. In particular, a subsymmetric basic sequence is
tail equivalent. \vskip .4cm

 \begin{theo} \label{main2}
 Let $X$ be a Banach space with a basis and having the property that for
any infinite dimensional  closed subspace $Y$ of $X$ there exist a
  constant $C_Y$ and infinite dimensional  subspaces $\mathcal E=(E_i)_{i=1}^{\infty}$ of $Y$ such that
  all normalized sequences $(x_n)_n$ with $x_n\in E_n$ for all $n$ are $C_Y$-tail equivalent. Then  $X$ contains a basic sequence equivalent
to the unit vector basis of $l_p$ or $c_0$.
\end{theo}

\noindent {\bf Proof. }
By passing to a subspace and relabeling assume that $X$ has a basis and $M$ is its basis
constant. According to Lemma~\ref{mainLemma} let $Y$ be a block subspace of $X$ such that
for every sequence $(Y_i)_{i=1}^\infty$ of infinite dimensional block subspaces of $Y$
there exists a normalized sequence $(y_i)_{i=1}^n$ with $y_i \in Y_i$ for $1 \leq i \leq n$
and $(y_i)_{i=1}^n$ is $2$-equivalent to the unit vector basis of $\ell_p^n$. By our
hypothesis there exists a sequence $\mathcal{E}= (E_j)$ of infinite dimensional subspaces
of $Y$  such that all normalized sequences $(x_n)$ with $x_n \in E_n$ are $C$-tail
equivalent. By standard perturbation arguments we can assume that ${\mathcal E}$ consists
of block subspaces.
 We then build inductively (by successive applications of Lemma~\ref{mainLemma})
a  normalized block sequence $(z_i)_i$ with $z_i\in G_i$
 for all $i$ such that $\{z_1,z_2\}$ is $2$-equivalent to the unit
 vector basis of $l_p^2$, $\{z_3,z_4, z_5\}$ is $2$-equivalent to the unit
 vector basis of $l_p^3$, $\{z_6,z_7, z_8, z_9\}$ is $2$-equivalent to the unit
 vector basis of $l_p^4$ and so on. Clearly $(z_i)_i$ is a block basic sequence
 with basis constant at most $M$.
 Fix $n\in\mathbb{N}$. From the definition of tail equivalence we
 can find $N$ large enough such that $(z_i)_{i=1}^{\infty}$ is
 $C$-equivalent to $(z_i)_{i=N+1}^{\infty}$ and $\{z_{N+1}, z_{N+2},\dots,
 z_{N+n}\}$ overlaps with at most two finite sequences of $z$'s as
 above. In other words, there exists $k\leq n$ such that
 $(z_i)_{y=N+1}^{N+k}$ is $2$-equivalent to the unit vector basis of
 $l_p^{k}$ and $(z_i)_{y=N+k+1}^{N+n}$ is $2$-equivalent to the unit vector basis of
 $l_p^{n-k}$. Then it easily
follows that $\{z_{N+1}, z_{N+2},\dots, z_{N+n}\}$ is $16(M+1)$-equivalent to the
unit vector basis of $l_p^{n}$. Since $(z_i)_{i=1}^{\infty}$ is
 $C$-equivalent to $(z_i)_{i=N+1}^{\infty}$, it follows that
$(z_i)_{i=1}^{n}$ is $16C(M+1)$-equivalent to the unit
vector basis of $l_p^{n}$, and this finishes the proof. \hfill $\blacksquare$

\footnotesize
\address


\begin{thebibliography}{A-A}
\frenchspacing


\bibitem[F-F-K-R]{FFKR} T. Figiel, R. Frankiewicz, R.A. Komorowski and C.
Ryll-Nardzewski, {\it Selecting basic sequences in $\vphi$-stable
Banach spaces}, Studia Math. {\bf 159 (3)} (2003), 499-515.

\bibitem[G] {G} W.T. Gowers, {\it An infinite Ramsey theorem and some Banach-space
dichotomies}, Ann. of Math. (2) {\bf 156} (2002), no.3, 797-833.


\bibitem[K]{K} J.L. Krivine, {\it Sous espaces de dimension finie
 des espaces de Banach r\'{e}ticul\'{e}s}, Ann. of Math. (2) {\bf
 104} (1976), 273-295.


\bibitem[M-TJ]{MT} V.D. Milman and N.
Tomczak-Jaegermann, {\it Asymptotic $l_p$ spaces and bounded
distorsions}, Contemp. Math. {\bf 144} (1993), 173-195.

\bibitem[Tc]{Tc} A. Tcaciuc {\it On the existence of asymptotic-$l_p$ structures in Banach spaces},
Canad. Math. Bull. {\bf 160(4)} (2007), 619-631

\bibitem[Z]{Z} M. Zippin, {\it On perfectly homogeneous bases in Banach
spaces}, Israel J. Math. {\bf 4} (1966), 265-272.


\end{thebibliography}
\end{document}